\documentclass[oneside, reqno]{amsart}

\usepackage{xypic}
\input xy
\xyoption{all}
\usepackage{epsfig}
\usepackage{amsthm}
\usepackage{amssymb}
\usepackage{amsmath}
\usepackage{amscd}
\usepackage{color}
\usepackage[T1]{fontenc}
\usepackage{stackengine}
\usepackage{upgreek}
\usepackage[font=scriptsize]{caption}
\usepackage{wrapfig}
\stackMath

\newcommand{\bg}{\begin{equation}}
\newcommand{\ed}{\end{equation}}
\newcommand{\bga}{\begin{eqnarray}}
\newcommand{\eda}{\end{eqnarray}}
\newcommand{\pf}{\textbf{Proof:\ }}

\def\cbdu{\par{\raggedleft$\Box$\par}}

\newtheorem {Theorem}  {Theorem}

\numberwithin{Theorem}{section}

\newtheorem {Lemma}[Theorem]  {Lemma}

\theoremstyle{definition}
\newtheorem{Definition}[Theorem]{Definition}
\theoremstyle{remark}
\newtheorem{Remark}[Theorem]{\bf Remark}

\newtheorem {Corollary}[Theorem]{\bf Corollary}

\expandafter\chardef\csname pre amssym.def
at\endcsname=\the\catcode`\@ \catcode`\@=11
\def\undefine#1{\let#1\undefined}
\def\newsymbol#1#2#3#4#5{\let\next@\relax
 \ifnum#2=\@ne\let\next@\msafam@\else
 \ifnum#2=\tw@\let\next@\msbfam@\fi\fi
 \mathchardef#1="#3\next@#4#5}
\def\mathhexbox@#1#2#3{\relax
 \ifmmode\mathpalette{}{\m@th\mathchar"#1#2#3}%
 \else\leavevmode\hbox{$\m@th\mathchar"#1#2#3$}\fi}
\def\hexnumber@#1{\ifcase#1 0\or 1\or 2\or 3\or 4\or 5\or 6\or 7\or 8\or
 9\or A\or B\or C\or D\or E\or F\fi}

\font\teneufm=eufm10 \font\seveneufm=eufm7 \font\fiveeufm=eufm5
\newfam\eufmfam
\textfont\eufmfam=\teneufm \scriptfont\eufmfam=\seveneufm
\scriptscriptfont\eufmfam=\fiveeufm

\catcode`\@=\csname pre amssym.def at\endcsname

\newcounter{remark}
\def  \12  {{\frac{1}{2}}}

\def\build#1_#2^#3{\mathrel{\mathop{\kern 0pt#1}\limits_{#2}^{#3}}}

\begin{document}

\title[Blow-up of dyadic MHD models]{Blow-up of dyadic MHD models with forward energy cascade}


\author [Mimi Dai]{Mimi Dai}

\address{Department of Mathematics, Statistics and Computer Science, University of Illinois at Chicago, Chicago, IL 60607, USA}
\email{mdai@uic.edu}

\thanks{The author was partially supported by NSF grants DMS--1815069 and DMS--2009422.}





\begin{abstract}
A particular type of dyadic model for the magnetohydrodynamics (MHD) with forward energy cascade is studied. The model includes intermittency dimension $\delta$ in the nonlinear scales. It is shown that when $\delta$ is small, positive solution with large initial data for either the dyadic MHD model or the dyadic Hall MHD develops blow-up in finite time.

\bigskip


\hspace{0.02cm}CLASSIFICATION CODE: 35Q35, 76D03, 76W05.
\end{abstract}

\maketitle

\section{Introduction}

Dyadic models for the incompressible magnetohydrodynamics (MHD) with Hall effect governed by 
\begin{equation}\label{hmhd}
\begin{split}
u_t+u\cdot\nabla u-B\cdot\nabla B+\nabla p=&\ \nu\Delta u,\\
B_t+u\cdot\nabla B-B\cdot\nabla u +d_i \nabla\times ((\nabla\times B)\times B)=&\ \mu\Delta B,\\
\nabla \cdot u=&\ 0,
\end{split}
\end{equation}
were derived in \cite{Dai-20}, where intermittency effect enters the derivation in a natural way. 
In system (\ref{hmhd}), the unknown functions $u$, $p$ and $B$ denote respectively the electrically conducting fluid velocity field, fluid pressure, and magnetic field influenced by the conducting fluid. 
The parameters $\nu,\mu$ and $d_i$ stand for the kinematic viscosity, magnetic resistivity and ion inertial length, respectively.  We assume (\ref{hmhd}) is posed either on $\mathbb R^3\times [0,\infty)$ or $\mathbb T^3\times [0,\infty)$. 
A general form of the derived dyadic (shell) model for (\ref{hmhd}) reads as
\begin{equation}\label{gen1}
\begin{split}
\frac{d}{dt}a_j+\nu\lambda_j^2 a_j\\
+\alpha_1\left(\lambda_j^{\frac{5-\delta_u}{2}}a_ja_{j+1}-\lambda_{j-1}^{\frac{5-\delta_u}{2}}a_{j-1}^2\right)+\beta_1\left(\lambda_{j}^{\frac{5-\delta_u}{2}}a_{j+1}^2-\lambda_{j-1}^{\frac{5-\delta_u}{2}}a_{j-1}a_j\right)\\
+\alpha_3\left(\lambda_{j}^{\frac{5-\delta_b}{2}}b_jb_{j+1}-\lambda_{j-1}^{\frac{5-\delta_b}{2}}b_{j-1}^2\right)
+\beta_3\left(\lambda_{j+1}^{\frac{5-\delta_b}{2}}b_{j+1}^2-\lambda_{j}^{\frac{5-\delta_b}{2}}b_{j-1}b_j\right)=0,
\end{split}
\end{equation}
\begin{equation}\label{gen2}
\begin{split}
\frac{d}{dt}b_j+\mu\lambda_j^2 b_j\\
+\alpha_2\left(\lambda_j^{\frac{5-\delta_b}{2}}a_jb_{j+1}-\lambda_{j-1}^{\frac{5-\delta_b}{2}}a_{j-1}b_{j-1}\right)
+\beta_2\left(\lambda_{j+1}^{\frac{5-\delta_b}{2}}a_{j+1}b_{j+1}-\lambda_{j}^{\frac{5-\delta_b}{2}}a_{j}b_{j-1}\right)\\
+\alpha_3\left(\lambda_{j}^{\frac{5-\delta_b}{2}}b_ja_{j+1}-\lambda_{j-1}^{\frac{5-\delta_b}{2}}a_{j-1}b_{j-1}\right)+\beta_3\left(\lambda_{j+1}^{\frac{5-\delta_b}{2}}b_{j+1}a_{j+1}-\lambda_{j}^{\frac{5-\delta_b}{2}}a_{j-1}b_{j}\right)\\
+d_i\alpha_4\left(\lambda_j^{\frac{7-\delta_b}{2}}b_jb_{j+1}-\lambda_{j-1}^{\frac{7-\delta_b}{2}}b_{j-1}^2\right)
+d_i\beta_4\left(\lambda_j^{\frac{7-\delta_b}{2}}b_{j+1}^2-\lambda_{j-1}^{\frac{7-\delta_b}{2}}b_{j}b_{j-1}\right)=0,
\end{split}
\end{equation}
for $j\geq 1$, which is an ODE system of infinitely many equations. In system (\ref{gen1})-(\ref{gen2}), the unknown functions $a_j$ and $b_j$ appear to be the kinetic energy and magnetic energy in the $j$-th shell, respectively, in the derivation. However, they can also be treated as Fourier coefficients of $u$ and $B$, respectively. By convention, we take $a_0=b_0=0$. The parameter $\lambda_j=\lambda^j$ stands for the wavenumber of the $j$-th shell for some $\lambda>1$. The parameters $\delta_u$ and $\delta_b$ represent intermittency dimension for the velocity field $u$ and magnetic field $B$, respectively, which are defined through the saturation level of Bernstein's inequality, see \cite{CD-Kol, Dai-20}. To be physically relevant, $\delta_u$ and $\delta_b$ take values in $[0,3]$. The situation of $\delta_u=\delta_b=3$ corresponds to the Kolmogorov regime, in which case both of the conducting flow and magnetic field flow are homogeneous, isotropic and self-similar. In the case of $\delta_u=\delta_b=0$, both flows are extremely inhomogeneous and singular. The parameters $\alpha_k$ and $\beta_k$ for $1\leq k\leq 4$ play essential roles in interpreting energy transfer among shells and the coupling relationship between the velocity field and magnetic field. They will be further discussed at a later time.

The dyadic model (\ref{gen1})-(\ref{gen2}) is derived under the following principles: (i) kinetic energy and magnetic energy are balanced through each shell; (ii) the total energy is conserved when $\nu=\mu=0$; (iii) only local interactions among shells are taken into account (in fact, only interactions with the first neighbor shells are employed here).  One can check that the total energy 
\begin{equation}\label{eq-total-energy}
E(t)=\frac12\sum_{j\geq 1}\left(a_j^2(t)+b_j^2(t)\right)
\end{equation}
is indeed formally conserved for the model with $\nu=\mu=0$ and any parameters $\alpha_k$ and $\beta_k$, $1\leq k\leq 4$. Moreover, the total energy is also formally conserved for the system with: (i) $\alpha_k=0$ and $\beta_k\neq 0$ for $1\leq k\leq 4$, in which case the dyadic model is the Obukov type;  
 (ii) $\beta_k=0$ and $\alpha_k\neq 0$ for $1\leq k\leq 4$, in which case the dyadic model is the Katz-Pavl\'ovic (KP) type, see \cite{KP, KZ}. 
It is important to notice that
the sign of the parameters $\alpha_k$ and $\beta_k$ determines the direction of energy transfer: positive sign indicates forward energy cascade, while negative sign indicates backward energy cascade.

Dyadic models for hydrodynamics governed by the Navier-Stokes equation (NSE) and Euler equation have been extensively studied, for instance, see \cite{BFM, BM, BMR, Bif, Ch, CF, CFP1, CFP2, CLT, DS, FP, Fri, Gle, JL, KP, KZ, LPPPV, Ob, OY, Wal}. (It would be quite challenging to list the complete literature on this topic; thus, the author gives up such attempt here.) They serve as approximating models for the true fluid flows, which reflect some most essential features of the turbulent flows. In fact, taking $b_j=0$ for $j\geq0$, $\alpha_1=1$ and $\beta_1=0$ in (\ref{gen1}), the model reduces to the KP dyadic model; instead, taking $b_j=0$ for $j\geq0$, $\alpha_1=0$ and $\beta_1=1$ in (\ref{gen1}) makes it to be the Obukov model. 
One major shortage of these dyadic models is that spatial complexity and geometry structures of the original flows are over simplified. Nevertheless, the study of dyadic models has provided important insights in the understanding of hydrodynamic turbulence. 

Dyadic models for the MHD turbulence were also introduced and studied by physicists, see \cite{Bis94, GLPG}, the recent article \cite{PSF} and references therein. The dyadic model (\ref{gen1})-(\ref{gen2}), derived based on harmonic analysis techniques and with intermittency effect included automatically, recovers some models from the physics community which will be pointed out below at the proper place. The main aim of proposing model (\ref{gen1})-(\ref{gen2}) is two-fold: (i) understand how the behaviours of solutions depend on the intermittency effect; (ii) explore how  different energy cascade and coupling relationships affect the dynamics.

In \cite{Dai-20}, the questions of well-posedness and finite-time blow-up were addressed for a special case of the model (\ref{gen1})-(\ref{gen2}). The particular model is obtained by taking the parameters $\alpha_1=\alpha_2=\alpha_4=1$, $\alpha_3=-1$, $\beta_k=0$ with $1\leq k\leq 4$, and $\delta_u=\delta_b:=\delta$ in (\ref{gen1})-(\ref{gen2}). Namely, the following model was considered,
\begin{equation}\label{sys-1}
\begin{split}
\frac{d}{dt}a_j+\nu\lambda_j^2 a_j=&\ -\left(\lambda_j^{\theta}a_ja_{j+1}-\lambda_{j-1}^{\theta}a_{j-1}^2\right)
+\left(\lambda_{j}^{\theta}b_jb_{j+1}-\lambda_{j-1}^{\theta}b_{j-1}^2\right) ,\\
\frac{d}{dt}b_j+\mu\lambda_j^2 b_j= &\ -\left(\lambda_j^{\theta}a_jb_{j+1}-\lambda_{j}^{\theta}b_ja_{j+1}\right)
-d_i\left(\lambda_j^{\theta+1}b_jb_{j+1}-\lambda_{j-1}^{\theta+1}b_{j-1}^2\right),
\end{split}
\end{equation}
for $j\geq 1$, $a_0=b_0=0$, and $\theta=\frac{5-\delta}{2}$. Some important features about this model are described below. First, the total energy $E(t)$ as in (\ref{eq-total-energy}) is formally conserved in the inviscid non-resistive case, i.e. $\nu=\mu=0$; the cross helicity defined by
\begin{equation}\label{eq-cross}
H^c(t)= \sum_{j\geq 1} a_j(t)b_j(t)
\end{equation}
is also formally conserved in the inviscid non-resistive MHD case, i.e. $\nu=\mu=d_i=0$. Second, from the signs of the nonlinear terms, we observe that energy moves towards larger frequency (forward energy cascade) within the nonlinear structures of $(u\cdot\nabla) u$ and the Hall effect $\nabla\times ((\nabla\times B)\times B)$; energy moves toward smaller frequency (backward energy cascade) within the coupled nonlinear structures. The following diagram illustrate the energy transfer among neighbor shells for system (\ref{sys-1}), 
\[ \cdot\cdot\cdot  \longrightarrow  \ \ a_{j-1} \ \ \longrightarrow \ \ a_j \ \ \longrightarrow \ \ a_{j+1} \ \ \longrightarrow \ \cdot\cdot\cdot\]
\[\hspace{-4mm}\uparrow \ \ \ \ \ \ \swarrow \ \  \ \ \uparrow \ \ \hspace{4mm} \swarrow \ \  \ \ \uparrow\]
\[ \cdot\cdot\cdot  \longrightarrow  \ \ b_{j-1} \ \ \longrightarrow \ \ b_j \ \ \longrightarrow \ \ b_{j+1} \ \ \longrightarrow \ \cdot\cdot\cdot\]
Third, when $\delta=3$ and hence $\theta=1$, this model corresponds to the so called L1 model derived by physicists, see \cite{GLPG, PSF}. The name L1 means that each flux term has local, two feet in the same shell and the third foot in a neighboring shell, for instance, $\lambda_j^\theta a_j^2a_{j+1}$ and $\lambda_j^\theta a_j b_jb_{j+1}$. 

For system (\ref{sys-1}), existence of global in time weak solutions is obtained in \cite{Dai-20} for any $\delta\in [0,3]$ (and any $\theta>0$);  when $d_i>0$, strong solution is shown to exist locally for $\delta\in (1, 3]$ and globally for $\delta= 3$; while for $d_i=0$, strong solution can be obtained locally for $\delta\in [0,3]$ and globally for $\delta\in [1,3]$. Moreover, when $d_i>0$ and $\delta<-1$, positive solution of (\ref{sys-1}) with large initial data is shown to develop blow-up at finite time. However, the last scenario is physically irrelevant since the intermittency dimension $\delta$ is only physically meaningful if $\delta\in [0,3]$. We point out that the question of whether positive solution of (\ref{sys-1}) with $d_i=0$ (i.e. the MHD dyadic model) develops blow-up at finite time remains open. 


In the current paper, we will work with another particular case of the general dyadic model (\ref{gen1})-(\ref{gen2}) with only forward energy cascade. Specifically, we assume $\delta_u=\delta_b=\delta$ and denote $\theta=\frac{5-\delta}{2}$ as before, take $\alpha_1=\alpha_3=\alpha_4=1$, $\alpha_2=-1$, $\beta_k=0$ with $1\leq k\leq 4$, and consider the following model
\begin{equation}\label{sys-2}
\begin{split}
\frac{d}{dt}a_j+\nu\lambda_j^2 a_j=&\ -\left(\lambda_j^{\theta}a_ja_{j+1}-\lambda_{j-1}^{\theta}a_{j-1}^2\right)
-\left(\lambda_{j}^{\theta}b_jb_{j+1}-\lambda_{j-1}^{\theta}b_{j-1}^2\right) ,\\
\frac{d}{dt}b_j+\mu\lambda_j^2 b_j= &\ \left(\lambda_j^{\theta}a_jb_{j+1}-\lambda_{j}^{\theta}b_ja_{j+1}\right)
-d_i\left(\lambda_j^{\theta+1}b_jb_{j+1}-\lambda_{j-1}^{\theta+1}b_{j-1}^2\right),
\end{split}
\end{equation}
for $j\geq 1$ and $a_0=b_0=0$. An obvious difference between system (\ref{sys-1}) and system (\ref{sys-2}) is the sign of the coupling terms $\left(\lambda_{j}^{\theta}b_jb_{j+1}-\lambda_{j-1}^{\theta}b_{j-1}^2\right)$ and $\left(\lambda_j^{\theta}a_jb_{j+1}-\lambda_{j}^{\theta}b_ja_{j+1}\right)$. That leads to some more sophisticated different features. For system (\ref{sys-2}), although the total energy is still formally conserved if $\nu=\mu=0$, the cross helicity as defined in (\ref{eq-cross}) is no longer conserved with $\nu=\mu=d_i=0$. Another important feature is that there is only forward energy cascade within the dynamics, see the illustration below
\[ \cdot\cdot\cdot  \longrightarrow  \ \ a_{j-1} \ \ \longrightarrow \ \ a_j \ \ \longrightarrow \ \ a_{j+1} \ \ \longrightarrow \ \cdot\cdot\cdot\]
\[\hspace{-4mm}\downarrow \ \ \ \ \ \ \nearrow \ \  \ \ \downarrow \ \ \hspace{4mm} \nearrow \ \  \ \ \downarrow\]
\[ \cdot\cdot\cdot  \longrightarrow  \ \ b_{j-1} \ \ \longrightarrow \ \ b_j \ \ \longrightarrow \ \ b_{j+1} \ \ \longrightarrow \ \cdot\cdot\cdot\]

The existence of short time strong solution and global strong solution to (\ref{sys-2}) can be established for proper regimes of the intermittency dimension $\delta$, in a similar fashion as the analysis for (\ref{sys-1}) in \cite{Dai-20}. In this paper,  we pursue to construct finite-time blow-up solutions to (\ref{sys-2}) with either $d_i=0$ or $d_i>0$ when the intermittency dimension is below certain threshold.  

In the case of the MHD dyadic model with $\theta>3$, we will show that blow-up develops at finite time for positive solution with large initial data in the space $H^s\times H^s$ with $s>\frac13\theta$.

\begin{Theorem}\label{thm-blow-mhd}
Let $(a(t),b(t))$ be a positive solution to (\ref{sys-2}) with $d_i=0$ and $\theta>3$. Let $\lambda\geq 2$. For any $\gamma>0$, there exists a constant $M_0$ such that if $\|a(0)\|_\gamma^2+\|b(0)\|_\gamma^2>M_0^2$, then $\|a(t)\|_{\frac13\theta+\frac23\gamma}^3+\|b(t)\|_{\frac13\theta+\frac23\gamma}^3$ is not locally integrable on $[0,\infty)$. 
\end{Theorem}

On the other hand, for the Hall MHD dyadic model with $\theta>3$, finite time blow-up occurs for positive solution with large initial data in the space $H^s\times H^{\frac13+s}$ with $s>\frac13\theta$.

\begin{Theorem}\label{thm-blow-hmhd}
Let $(a(t),b(t))$ be a positive solution to (\ref{sys-2}) with $d_i>0$ and $\theta>3$. For any $\gamma>0$, there exists a constant $M_0$ such that if $\|a(0)\|_\gamma^2+\|b(0)\|_\gamma^2>M_0^2$, then $\|a(t)\|_{\frac13\theta+\frac23\gamma}^3+\|b(t)\|_{\frac13(\theta+1)+\frac23\gamma}^3$ is not locally integrable on $[0,\infty)$.
\end{Theorem}

\begin{Remark}\label{rk1}
Since $\theta=\frac{5-\delta}{2}$, $\theta>3$ is equivalent to $\delta<-1$. 
\end{Remark}

\begin{Remark}\label{rk2}
In Theorem \ref{thm-blow-mhd}, the parameter $\lambda$, the basis of the wavenumber $\lambda_j=\lambda^j$, can be taken as any value larger than 1.  To reduce the complexity of analyzing parameters satisfying (\ref{para-11})-(\ref{para-16}), we choose $\lambda\geq 2$. 
\end{Remark}

\begin{Remark}\label{rk3}
The question whether a solution of (\ref{sys-2}) with positive initial data remains positive is open and will be addressed in future investigation. It is known that, if $B=0$ and hence $b_j=0$ for all $j\geq 0$, the reduced NSE dyadic model (\ref{sys-2}) with positive initial data produces positive solutions, see \cite{Ch}. 
\end{Remark}

\begin{Remark}\label{rk4}
In view of the fact that the Hall MHD system (\ref{hmhd}) with $d_i>0$ involves a more singular nonlinear structure of the Hall effect,  reflected in the dyadic model (\ref{gen1})-(\ref{gen2}) with a larger nonlinear scale $d_i\left(\lambda_j^{\theta+1}b_jb_{j+1}-\lambda_{j-1}^{\theta+1}b_{j-1}^2\right)$, one might expect to show blow-up for system (\ref{sys-2}) with $d_i>0$ for smaller $\theta$, that is, for $\theta<3$. Nevertheless, in the proof of Theorem \ref{thm-blow-hmhd} in Section \ref{sec-blow-hmhd}, it appears that the coupling terms cause serious barrier to lower the threshold of $\theta$ for blow-up. That could be just the limitation of the approach of proving blow-up in this paper. There is hope to move down the threshold of $\theta$ for blow-up by other frameworks of proving blow-up.
\end{Remark}



An interesting connection between the intermittency effect and dissipation strength can be revealed through dyadic models in the following way. In fact, the dyadic system (\ref{sys-2}) can be rescaled to
\begin{equation}\label{sys-3}
\begin{split}
\frac{d}{dt}a_j+\nu\bar\lambda_j^{2\alpha} a_j=& -\bar\lambda_j a_ja_{j+1}+\bar\lambda_{j-1} a_{j-1}^2-\bar\lambda_j b_jb_{j+1}+\bar\lambda_{j-1} b_{j-1}^2,\\
\frac{d}{dt}b_j +\mu\bar\lambda_j^{2\alpha} b_j=&\ \bar\lambda_j a_jb_{j+1}-\bar\lambda_j b_ja_{j+1}-d_i\left(\bar\lambda_j^{\alpha+1}b_jb_{j+1}-\bar\lambda_{j-1}^{\alpha+1}b_{j-1}^2\right)
\end{split}
\end{equation}
with 
\[\alpha=\frac1{\theta}=\frac{2}{5-\delta},\] 
by rescaling the wavenumber $\lambda_j=\bar\lambda_j^{\alpha}$. System (\ref{sys-3}) can be seen as the dyadic model of the Hall-MHD system with generalized diffusions $(-\Delta)^\alpha u$ and $(-\Delta)^\alpha B$. 
The results of Theorem \ref{thm-blow-mhd} and Theorem \ref{thm-blow-hmhd} can be transformed to system (\ref{sys-3}) as follows.

\begin{Corollary}\label{cor-blow-mhd}
Let $(a(t),b(t))$ be a positive solution to (\ref{sys-3}) with $d_i=0$ and $\alpha<\frac13$. 
For any $\gamma>0$, there exists a constant $M_0$ such that if $\|a(0)\|_\gamma^2+\|b(0)\|_\gamma^2>M_0^2$, then $\|a(t)\|_{\frac13+\gamma}^3+\|b(t)\|_{\frac13+\gamma}^3$ is not locally integrable on $[0,\infty)$. 
\end{Corollary}

\begin{Corollary}\label{cor-blow-hmhd}
Let $(a(t),b(t))$ be a positive solution to (\ref{sys-3}) with $d_i>0$ and $\alpha<\frac13$. For any $\gamma>0$, there exists a constant $M_0$ such that if $\|a(0)\|_\gamma^2+\|b(0)\|_\gamma^2>M_0^2$, then $\|a(t)\|_{\frac13+\gamma}^3+\|b(t)\|_{\frac13(\alpha+1)+\gamma}^3$ is not locally integrable on $[0,\infty)$.
\end{Corollary}

The proof of Theorem \ref{thm-blow-mhd} and Theorem \ref{thm-blow-hmhd} relies on a contradiction argument and the construction of a Lyapunov function $\mathcal L(t)$ which would satisfy a Riccati type of inequality. Depending on whether $d_i>0$ or not, i.e. whether the Hall term is present or not, the choice of $\mathcal L(t)$ is different. The construction of $\mathcal L(t)$ for both the dyadic MHD and Hall MHD models is described in Section \ref{sec-est}; some properties of $\mathcal L(t)$ are also established there. The proof of Theorem \ref{thm-blow-mhd} and Theorem \ref{thm-blow-hmhd} is provided in Section \ref{sec-blow-mhd} and Section \ref{sec-blow-hmhd}, respectively. On the other hand, Corollary \ref{cor-blow-mhd} and Corollary \ref{cor-blow-hmhd} can be justified automatically from the rescaling relationship.

\section{Notations and notion of solutions}
We denote $H=l^2$ which is endowed with the standard scalar product and norm,
\[(u,v):=\sum_{n=1}^\infty u_nv_n, \ \ \ |u|:=\sqrt{(u,u)}.\]
As mentioned earlier, we choose the wavenumber $\lambda_n=\lambda^n$ for a constant $\lambda>1$, and all integers $n\geq 1$.
Corresponding to the standard Sobolev space $H^s$ for functions with spacial variables, we
use the same notation $H^s$ here to represent the space for a sequence $\{u_n\}_{n=1}^{\infty}$, which is endowed with the scaler product 
\[(u,v)_s:=\sum_{n=1}^\infty \lambda_n^{2s}u_nv_n\]
and the norm
\[\|u\|_s:=\sqrt{(u,u)_s}.\]
We notice that $H^0=H=l^2$ which is regarded as the energy space. 

In the following, the concept of solutions for the dyadic system (\ref{sys-2}) is introduced. 
\begin{Definition}\label{def1}
A pair of $H$-valued functions $(a(t), b(t))$ defined on $[t_0,\infty)$ is said to be a weak solution of (\ref{sys-2}) if $a_j$ and $b_j$ satisfy (\ref{sys-2}) and $a_j, b_j\in C^1([t_0,\infty))$ for all $j\geq0$.
\end{Definition}

\begin{Definition}\label{def2}
A solution $(a(t), b(t))$ of (\ref{sys-2}) is strong on $[T_1, T_2]$ if $\|a\|_1$ and $\|b\|_1$ are bounded on $[T_1, T_2]$. A solution is strong on $[T_1, \infty)$ if it is strong on every interval $[T_1, T_2]$ for any $T_2>T_1$.
\end{Definition}


\section{Lyapunov function and auxiliary estimates}
\label{sec-est}

In this section, we construct a Lyapunov function for system (\ref{sys-2}) and present its continuity under certain conditions. In particular, if $d_i=0$, we consider 
\begin{equation}\label{L1}
\begin{split}
\mathcal L(t):= &\ \|a(t)\|_\gamma^2+\|b(t)\|_\gamma^2+c_1\sum_{j=1}^\infty \lambda_j^{2\gamma}a_j(t)a_{j+1}(t)\\
&+c_2\sum_{j=1}^\infty \lambda_j^{2\gamma}b_j(t)a_{j+1}(t)+c_3\sum_{j=1}^\infty \lambda_j^{2\gamma}a_j(t)b_j(t)
\end{split}
\end{equation}
for some appropriate positive constants $c_1, c_2$, and $c_3$.  
The main principle of designing $\mathcal L(t)$ is to have terms $a_j^3$ and $b_j^3$ included in $\frac{d}{dt} \mathcal L(t)$, which will play a crucial role to derive a Riccati type of inequality for $\mathcal L(t)$. In fact, $\frac{d}{dt} (a_ja_{j+1})$ produces $\lambda_j^{\theta}a_j^3$ and $\frac{d}{dt} (b_ja_{j+1})$ gives $\lambda_j^{\theta}b_j^3$. However, it turns out that the term $\lambda_j^{\theta}b_j^3$ is not enough to control a flux triple term $\lambda_j^{\theta}b_jb_{j+1}b_{j+2}$ in the estimates. It is the reason that we include the term $\lambda_j^{2\gamma}a_jb_j$ in $\mathcal L(t)$, and hence $\frac{d}{dt} (a_jb_j)$ gives a term $\lambda_j^{\theta}b_j^2b_{j+1}$ which can contribute to control $\lambda_j^{\theta}b_jb_{j+1}b_{j+2}$.

For the dyadic Hall MHD model (\ref{sys-2}) with $d_i>0$, we choose
\begin{equation}\label{L2}
\begin{split}
\mathcal L(t):= &\ \|a(t)\|_\gamma^2+\|b(t)\|_\gamma^2+c_1\sum_{j=1}^\infty \lambda_j^{2\gamma}a_j(t)a_{j+1}(t)\\
&+c_2\sum_{j=1}^\infty \lambda_j^{2\gamma}b_j(t)b_{j+1}(t)
\end{split}
\end{equation}
for appropriate constants $c_1>0$ and $c_2>0$.  As for the MHD case, $\frac{d}{dt} (a_ja_{j+1})$ includes the good term $\lambda_j^{\theta}a_j^3$. While in this case, $\frac{d}{dt} (b_jb_{j+1})$ contributes a good term $\lambda_j^{\theta+1}b_j^3$ due to the presence of the Hall term; in the same time, $\frac{d}{dt} b_j^2$ (from $\frac{d}{dt} \|b(t)\|_\gamma^2$) gives  $\lambda_j^{\theta+1}b_j^2b_{j+1}$, also due to the Hall effect. Thus, the two good terms together are able to control many negative flux terms including $\lambda_j^{\theta+1}b_jb_{j+1}b_{j+2}$.

In the rest of this section, we will provide some auxiliary estimates and show the continuity of $\mathcal L(t)$ under certain conditions on the solution.

\begin{Lemma}\label{le-triple}
(i) If $\theta>3+\gamma$, there exists a constant $c_0>0$ such that
\[\sum_{j=1}^\infty \lambda_j^{2\gamma+\theta}a_j^3\geq c_0\|a\|_{\gamma+1}^3, \ \ \sum_{j=1}^\infty \lambda_j^{2\gamma+\theta}b_j^3\geq c_0\|b\|_{\gamma+1}^3, \]
\[\sum_{j=1}^\infty \lambda_j^{2\gamma+\theta+1}b_j^3\geq c_0\|b\|_{\gamma+1}^3.\]
(ii) If $\theta>3+\gamma$, we also have
\[\|a(t)\|_{\gamma+1}\leq \|a(t)\|_{\frac13\theta+\frac23\gamma}, \ \ \|b(t)\|_{\gamma+1}\leq \|b(t)\|_{\frac13\theta+\frac23\gamma}.\] 
(iii) The following inequalities 
\[\sum_{j=1}^\infty \lambda_j^{2\gamma+2}a_ja_{j+1}\leq \lambda^{-\gamma-1}\|a\|_{\gamma+1}^2\] 
\[\sum_{j=1}^\infty \lambda_j^{2\gamma+2}b_jb_{j+1}\leq \lambda^{-\gamma-1}\|b\|_{\gamma+1}^2\]
\[\sum_{j=1}^\infty \lambda_j^{2\gamma+2}a_jb_j\leq \frac12 \left(\|a\|_{\gamma+1}^2+\|b\|_{\gamma+1}^2\right)\]
\[\sum_{j=1}^\infty \lambda_j^{2\gamma+2}b_ja_{j+1}\leq \frac12 \lambda^{-\gamma-1}\left(\|a\|_{\gamma+1}^2+\|b\|_{\gamma+1}^2\right)\]
 hold. \\
\noindent (iv) For positive $a_j$ and $b_j$ with $j\geq 1$, we have
 \[\sum_{j=1}^\infty\lambda_j^{2\gamma+\theta}a_j^2a_{j+1}\leq 2 \|a\|_{\frac13\theta+\frac23\gamma}^3,\]
 \[\sum_{j=1}^\infty\lambda_j^{2\gamma+\theta}b_j^2b_{j+1}\leq 2 \|b\|_{\frac13\theta+\frac23\gamma}^3,\]
 \[\sum_{j=1}^\infty\lambda_j^{2\gamma+\theta}b_j^2a_{j+1}\leq \|a\|_{\frac13\theta+\frac23\gamma}^3+ \|b\|_{\frac13\theta+\frac23\gamma}^3.\]

\end{Lemma}
\pf
The justification of the inequalities in (i) is rather standard and thus omitted here. One can find a quick proof in \cite{Dai-20}. The inequalities in (ii) follow immediately from the fact $\theta>3+\gamma$ and hence $\gamma+1<\frac13\theta+\frac23\gamma$. The ones in (iii) are not complicated either and we only show one of them below. Applying H\"older's and Young's inequalities, we have
\begin{equation}\notag
\begin{split}
\sum_{j=1}^\infty \lambda_j^{2\gamma+2}b_ja_{j+1}=&\lambda^{-\gamma-1}\sum_{j=1}^\infty\left(\lambda_j^{\gamma+1}b_j\right) \left(\lambda_{j+1}^{\gamma+1}a_{j+1}\right)\\
\leq& \lambda^{-\gamma-1}\left(\sum_{j=1}^\infty\lambda_j^{2\gamma+2}b_j^2\right)^{\frac12}
\left(\sum_{j=1}^\infty\lambda_{j+1}^{2\gamma+2}a_{j+1}^2\right)^{\frac12}\\
\leq&\frac12 \lambda^{-\gamma-1}\left(\|a\|_{\gamma+1}^2+\|b\|_{\gamma+1}^2\right).
\end{split}
\end{equation}

We only show the last inequality of (iv); another two can be proved similarly. The application of Young's inequality and a basic inequality for sum leads to
\begin{equation}\notag
\begin{split}
\sum_{j=1}^\infty\lambda_j^{2\gamma+\theta}b_j^2a_{j+1}\leq &\sum_{j=1}^\infty\lambda_j^{2\gamma+\theta}\left(\frac23b_j^3+\frac13a_{j+1}^3\right)\\
\leq & \sum_{j=1}^\infty\lambda_j^{2\gamma+\theta}a_j^3+\sum_{j=1}^\infty\lambda_j^{2\gamma+\theta}b_j^3\\
\leq & \left(\sum_{j=1}^\infty\lambda_j^{\frac23(2\gamma+\theta)}a_j^2\right)^{\frac32}+\left(\sum_{j=1}^\infty\lambda_j^{\frac23(2\gamma+\theta)}b_j^2\right)^{\frac32}\\
\leq &\|a\|_{\frac13\theta+\frac23\gamma}^3+\|b\|_{\frac13\theta+\frac23\gamma}^3.
\end{split}
\end{equation}

\cbdu

\begin{Lemma}\label{le-cont1}
 Let $(a(t),b(t))$ be a positive solution to (\ref{sys-2}) with $d_i=0$.
 Assume $\|a(t)\|_{\frac13\theta+\frac23\gamma}^3+\|b(t)\|_{\frac13\theta+\frac23\gamma}^3$ is locally integrable on $[0,\infty)$.
 Then $\mathcal L(t)$ defined in (\ref{L1}) is continuous on $[0,\infty)$. 
\end{Lemma}
\pf
We denote
\[E_\gamma(t):=\|a(t)\|_{\gamma}^2+\|b(t)\|_{\gamma}^2,  \]
\[f(t):=c_1\sum_{j=1}^\infty \lambda_j^{2\gamma}a_j(t)a_{j+1}(t)+c_2\sum_{j=1}^\infty \lambda_j^{2\gamma}b_j(t)a_{j+1}(t)+c_3\sum_{j=1}^\infty \lambda_j^{2\gamma}a_j(t)b_j(t).\]
Under the assumption, we show that both $E_\gamma$ and $f$ are continuous on $[0,\infty)$.

Applying the two equations of (\ref{sys-2}) with $d_i=0$, and taking the sum for all $j\geq 1$, we find that
\begin{equation}\notag
\begin{split}
&E_\gamma(t)-E_\gamma(0)\\
=&-2\int_0^t\nu \|a(\tau)\|_{\gamma+1}^2+\mu\|b(\tau)\|_{\gamma+1}^2\,d\tau
+2(\lambda^{2\gamma}-1)\int_0^t\sum_{j=1}^\infty \lambda_j^{2\gamma+\theta}a_j^2a_{j+1}\, d\tau\\
&+2(\lambda^{2\gamma}-1)\int_0^t\sum_{j=1}^\infty \lambda_j^{2\gamma+\theta}b_j^2a_{j+1}\, d\tau.
\end{split}
\end{equation}
Combining the inequalities of Lemma \ref{le-triple} (ii) and (iv)
and the assumption that $\|a(t)\|_{\frac13\theta+\frac23\gamma}^3+\|b(t)\|_{\frac13\theta+\frac23\gamma}^3$ is locally integrable, we conclude that $\|a(t)\|_{\gamma+1}^2$, $\|b(t)\|_{\gamma+1}^2$,  $\sum_{j=1}^\infty\lambda_j^{2\gamma+\theta}a_j^2a_{j+1}$, and
$\sum_{j=1}^\infty\lambda_j^{2\gamma+\theta}b_j^2a_{j+1}$ are all locally integrable as well. Therefore, the integrals on the right hand side of the equation above are all defined for any $t>0$. It thus follows that $E_\gamma$ is continuous on $[0,\infty)$.

We denote for $j\geq 1$
\[f_j(t)=c_1\lambda_j^{2\gamma}a_j(t)a_{j+1}(t)+c_2\lambda_j^{2\gamma}b_j(t)a_{j+1}(t)+c_3\lambda_j^{2\gamma}a_j(t)b_j(t). \]
For any $t_0>0$, we infer
\begin{equation}\label{lim0}
\begin{split}
&\limsup_{t\to t_0}\left|f(t)-f(t_0)\right|\\
=& \limsup_{t\to t_0}\left|\sum_{j=1}^\infty f_j(t)-\sum_{j=1}^\infty f_j(t_0)\right|\\
=&\lim_{J\to \infty} \limsup_{t\to t_0}\left| \sum_{j=1}^{J-1} f_j(t)-\sum_{j=1}^{J-1}f_j(t_0) +\sum_{j=J}^\infty f_j(t)-\sum_{j=J}^\infty f_j(t_0)\right|\\
\leq & \lim_{J\to \infty} \limsup_{t\to t_0} \sum_{j=1}^{J-1}\left| f_j(t)-f_j(t_0) \right|  + \lim_{J\to \infty} \limsup_{t\to t_0} \left|\sum_{j=J}^\infty f_j(t)-\sum_{j=J}^\infty f_j(t_0)\right|.\\
\end{split}
\end{equation}
According to the definition of solution in Definition \ref{def1}, $f_j$ is continuous for any $j\geq 1$ and hence 
\[\lim_{t-t_0} \left| f_j(t)-f_j(t_0) \right| =0, \ \ \ \ \forall \ \ \ 1\leq j\leq J-1.\]
It implies that 
\begin{equation}\label{lim1}
\lim_{J\to \infty} \limsup_{t\to t_0} \sum_{j=1}^{J-1}\left| f_j(t)-f_j(t_0) \right|=0.
\end{equation}
To analyze the last limit in (\ref{lim0}), we observe that from Lemma \ref{le-triple} (iii) 
\[0\leq f(t)\leq 4c_1\|a(t)\|_\gamma^2 +4c_2\|b(t)\|_\gamma^2\leq 4(c_1+c_2)E_\gamma(t). \]
The continuity of $E_\gamma $ on $[0,\infty)$ implies $f$ is bounded on every interval $[T_1,T_2]$, for any $T_2>T_1\geq 0$. Therefore, it follows that
\begin{equation}\label{lim2}
 \lim_{J\to \infty} \limsup_{t\to t_0} \left|\sum_{j=J}^\infty f_j(t)-\sum_{j=J}^\infty f_j(t_0)\right|=0
\end{equation}
In view of (\ref{lim0})-(\ref{lim2}), we claim $f$ is continuous on $[0,\infty)$. It accomplishes the proof of the lemma.

\cbdu

When $d_i>0$, we have the following statement. 

\begin{Lemma}\label{le-cont2}
 Let $(a(t),b(t))$ be a positive solution to (\ref{sys-2}) with $d_i>0$.
 Assume $\|a(t)\|_{\frac13\theta+\frac23\gamma}^3+\|b(t)\|_{\frac13(\theta+1)+\frac23\gamma}^3$ is locally integrable on $[0,\infty)$.
 Then $\mathcal L(t)$ defined in (\ref{L2}) is continuous on $[0,\infty)$. 
\end{Lemma}
\pf
The proof follows a close line to that of Lemma \ref{le-cont1}. We only explain why it requires $\|b(t)\|_{\frac13(\theta+1)+\frac23\gamma}^3$ to be locally integrable on $[0,\infty)$. 
Indeed, multiplying the $a_j$ equation of (\ref{sys-2}) by $\lambda_j^{2\gamma}a_j$ and the $b_j$ equation with $d_i>0$ by $\lambda_j^{2\gamma}b_j$, adding all the shells for $j\geq 1$, and integrating from $0$ to $t$,  we obtain
\begin{equation}\notag
\begin{split}
&E_\gamma(t)-E_\gamma(0)\\
=&-2\int_0^t\nu \|a(\tau)\|_{\gamma+1}^2+\mu\|b(\tau)\|_{\gamma+1}^2\,d\tau
+2(\lambda^{2\gamma}-1)\int_0^t\sum_{j=1}^\infty \lambda_j^{2\gamma+\theta}a_j^2a_{j+1}\, d\tau\\
&+2d_i(\lambda^{2\gamma}-1)\int_0^t\sum_{j=1}^\infty \lambda_j^{2\gamma+\theta+1}b_j^2b_{j+1}\, d\tau
+2(\lambda^{2\gamma}-1)\int_0^t\sum_{j=1}^\infty \lambda_j^{2\gamma+\theta}b_j^2a_{j+1}\, d\tau.
\end{split}
\end{equation}
Referring to the second inequality of Lemma \ref{le-triple} (iv) with $\theta$ replaced by $\theta+1$, the assumption $\|b(t)\|_{\frac13(\theta+1)+\frac23\gamma}^3$ is locally integrable on $[0,\infty)$ guarantees that $\sum_{j=1}^\infty \lambda_j^{2\gamma+\theta+1}b_j^2b_{j+1}$ is locally integrable on $[0,\infty)$.

\cbdu

\section{Blow-up of positive solutions of dyadic MHD}
\label{sec-blow-mhd}


This section is devoted to a proof of Theorem \ref{thm-blow-mhd}. The following lemma plays an important role, whose proof is postponed to a later time. 

\begin{Lemma}\label{le-blow-mhd}
Consider system (\ref{sys-2}) with $d_i=0$. Let $\theta>3$ and $0<\gamma<\theta-3$. Fix $\lambda\geq2$.
Assume $\|a(0)\|_\gamma^2+\|b(0)\|_\gamma^2> M_0^2$ for a certain constant $M_0>0$. Then, the function $\mathcal L(t)$ defined in (\ref{L1}) for positive solution $(a(t), b(t))$ of (\ref{sys-2}) is a Lyapunov function and it blows up in finite time. 
\end{Lemma}



\textbf{Proof of Theorem \ref{thm-blow-mhd}: }
We adapt a contradiction argument here. Suppose that $(a(t),b(t))$ is a positive solution to (\ref{sys-2}) with $d_i=0$ such that $\|a(t)\|_{\frac13\theta+\frac23\gamma}^3+\|b(t)\|_{\frac13\theta+\frac23\gamma}^3$ is integrable on $[0,T]$ for any $T>0$, with $\gamma\in(0,\theta-3)$ and $\gamma\ll 1$. It follows from Lemma \ref{le-cont1} that $\mathcal L$ defined in (\ref{L1}) is continuous on $[0,\infty)$. The assumption of $\|a(0)\|_\gamma^2+\|b(0)\|_\gamma^2> M_0^2$ for a constant $M_0>0$ along with Lemma \ref{le-blow-mhd} implies that the function $\mathcal L$ blows up in finite time. Obviously, the last two properties of $\mathcal L$ leads to a contradiction. 

We also notice that $\|a(t)\|_{\frac13\theta+\frac23\gamma}^3+\|b(t)\|_{\frac13\theta+\frac23\gamma}^3$ is not locally integrable for an arbitrarily small $\gamma>0$ implies it is not locally integrable for any $\gamma>0$. 

\cbdu

We are left to give a justification of Lemma \ref{le-blow-mhd}.

\textbf{Proof of Lemma \ref{le-blow-mhd}: } 
The goal is to show that for some $T>0$,
\[\mathcal L(t)>\mathcal L(0), \ \ \forall t\in(0,T],\]
and $\mathcal L$ satisfies a Riccati type of inequality.
 
Utilizing the two equations of (\ref{sys-2}) with $d_i=0$, straightforward computation shows that 
\begin{equation}\label{aa-1}
\begin{split}
\frac{d}{dt}\left(\lambda_j^{2\gamma}a_ja_{j+1}\right)=&-\nu(1+\lambda^2)\lambda_j^{2\gamma+2}a_ja_{j+1} +\lambda_j^{2\gamma+\theta} a_j^3+\lambda_j^{2\gamma+\theta} a_jb_j^2\\
&+\lambda_{j-1}^\theta \lambda_j^{2\gamma} a_{j-1}^2a_{j+1}
+\lambda_{j-1}^\theta \lambda_j^{2\gamma}b_{j-1}^2a_{j+1}\\
&-\lambda_j^{2\gamma+\theta} a_ja_{j+1}^2
-\lambda_j^{2\gamma}\lambda_{j+1}^\theta a_ja_{j+1}a_{j+2} \\
&-\lambda_j^{2\gamma+\theta} b_j a_{j+1}b_{j+1}-\lambda_j^{2\gamma}\lambda_{j+1}^\theta a_jb_{j+1}b_{j+2},\\
\end{split}
\end{equation}
\begin{equation}\label{ab-1}
\begin{split}
\frac{d}{dt}\left(\lambda_j^{2\gamma}b_ja_{j+1}\right)=&-(\mu+\nu\lambda^2)\lambda_j^{2\gamma+2}b_ja_{j+1}+\lambda_j^{2\gamma+\theta}b_j^3+\lambda_j^{2\gamma+\theta}a_j^2b_j\\
&+\lambda_j^{2\gamma+\theta}a_ja_{j+1}b_{j+1}-\lambda_j^{2\gamma+\theta}b_ja_{j+1}^2\\
&-\lambda_j^{2\gamma}\lambda_{j+1}^\theta b_ja_{j+1}a_{j+2}-\lambda_j^{2\gamma}\lambda_{j+1}^\theta b_jb_{j+1}b_{j+2},
\end{split}
\end{equation}
\begin{equation}\label{ab-j}
\begin{split}
\frac{d}{dt}\left(\sum_{j\geq 1}\lambda_j^{2\gamma}a_jb_j\right)=&-(\nu+\mu)\sum_{j\geq 1}\lambda_j^{2\gamma+2}a_jb_j+\left(\lambda^{2\gamma}+1\right)\sum_{j\geq 1}\lambda_j^{2\gamma+\theta}a_j^2b_{j+1}\\
&+\left(\lambda^{2\gamma}-1\right)\sum_{j\geq 1}\lambda_j^{2\gamma+\theta}b_j^2b_{j+1}
-2\sum_{j\geq 1}\lambda_j^{2\gamma+\theta}a_jb_ja_{j+1}.
\end{split}
\end{equation}
In the same time, we have the energy equality
\begin{equation}\label{a2b2}
\begin{split}
&\frac{d}{dt}\left(\|a(t)\|_\gamma^2+\|b(t)\|_\gamma^2\right)\\
=&\ -2\nu\|a(t)\|_{\gamma+1}^2-2\mu \|b(t)\|_{\gamma+1}^2+2(\lambda^{2\gamma}-1)\sum_{j=1}^\infty \lambda_j^{2\gamma+\theta}a_j^2a_{j+1}\\
&+2(\lambda^{2\gamma}-1)\sum_{j=1}^\infty \lambda_j^{2\gamma+\theta} b_{j}^2a_{j+1}.\\
\end{split}
\end{equation}
The task is to control the negative terms on the right hand side of (\ref{aa-1})-(\ref{a2b2}) using the positive terms $\lambda_j^{2\gamma+\theta} a_j^3$, $\lambda_j^{2\gamma+\theta} a_jb_j^2$, $\lambda_j^{2\gamma+\theta} b_j^3$, $\lambda_j^{2\gamma+\theta} a_j^2a_{j+1}$, $\lambda_j^{2\gamma+\theta} b_j^2a_{j+1}$ and $\lambda_j^{2\gamma+\theta} b_j^2b_{j+1}$.
We estimate these negative terms by applying Young's inequality as follows,
\begin{equation}\label{basic-ineq1}
\begin{split}
&\lambda_j^{2\gamma+\theta}a_ja_{j+1}^2\\
=&\ \lambda^{-\frac12(2\gamma+\theta)}\left(\lambda_j^{\frac12(2\gamma+\theta)}a_ja_{j+1}^{\frac12}\right)\left(\lambda_{j+1}^{\frac12(2\gamma+\theta)}a_{j+1}^{\frac32}\right)\\
\leq &\frac12\lambda^{-\frac12(2\gamma+\theta)} \lambda_{j+1}^{2\gamma+\theta}a_{j+1}^3+\frac12\lambda^{-\frac12(2\gamma+\theta)} \lambda_{j}^{2\gamma+\theta}a_j^2a_{j+1};
\end{split}
\end{equation}
\begin{equation}\label{basic-ineq2}
\begin{split}
&\lambda_j^{2\gamma}\lambda_{j+1}^\theta a_ja_{j+1}a_{j+2}\\
=&\ \lambda_j^{2\gamma}\lambda_{j+1}^\theta \left(a_ja_{j+1}^{\frac12}\right)  \left(a_{j+1}^{\frac12}a_{j+2}^{\frac14}\right)  \left(a_{j+2}^{\frac34}\right)\\
\leq &\ \frac12\lambda^{\theta}\lambda_j^{2\gamma+\theta} a_j^2a_{j+1}
+\frac14 \lambda^{-2\gamma}\lambda_{j+1}^{2\gamma+\theta} a_{j+1}^2a_{j+2}+\frac14 \lambda^{-4\gamma-\theta}\lambda_{j+2}^{2\gamma+\theta} a_{j+2}^3;
\end{split}
\end{equation}
\begin{equation}\label{basic-ineq3}
\begin{split}
&\lambda_j^{2\gamma+\theta} b_ja_{j+1}b_{j+1}\\
=&\ \lambda^{-\frac12(2\gamma+\theta)}\left(\lambda_j^{\frac12(2\gamma+\theta)}b_ja_{j+1}^{\frac12}\right)\left(\lambda_{j+1}^{\frac12(2\gamma+\theta)}a_{j+1}^{\frac12}b_{j+1}\right)\\
\leq&\ \frac12\lambda^{-\frac12(2\gamma+\theta)}\lambda_j^{2\gamma+\theta}b_j^2a_{j+1}+\frac12\lambda^{-\frac12(2\gamma+\theta)}\lambda_{j+1}^{2\gamma+\theta}a_{j+1}b_{j+1}^2;
\end{split}
\end{equation}
\begin{equation}\label{basic-ineq4-1}
\begin{split}
&\lambda_j^{2\gamma}\lambda_{j+1}^{\theta}a_jb_{j+1}b_{j+2}\\
=&\ \lambda^{-2\gamma}\left(\lambda_j^{\frac13(2\gamma+\theta)}a_j\right)\left(\lambda_{j+1}^{\frac13(2\gamma+\theta)}b_{j+1}\right)\left(\lambda_{j+2}^{\frac13(2\gamma+\theta)}b_{j+2}\right)\\
\leq& \frac13 \lambda^{-2\gamma}\lambda_j^{2\gamma+\theta}a_j^3+\frac13 \lambda^{-2\gamma}\lambda_{j+1}^{2\gamma+\theta}b_{j+1}^3+\frac13 \lambda^{-2\gamma}\lambda_{j+2}^{2\gamma+\theta}b_{j+2}^3;
\end{split}
\end{equation}
\begin{equation}\label{basic-ineq5}
\begin{split}
&\lambda_j^{2\gamma+\theta} b_ja_{j+1}^2\\
=&\ \lambda^{-\frac12(2\gamma+\theta)}\left(\lambda_j^{\frac12(2\gamma+\theta)}b_ja_{j+1}^{\frac12}\right)\left(\lambda_{j+1}^{\frac12(2\gamma+\theta)}a_{j+1}^{\frac32}\right)\\
\leq&\ \frac12\lambda^{-\frac12(2\gamma+\theta)}\lambda_j^{2\gamma+\theta}b_j^2a_{j+1}+\frac12\lambda^{-\frac12(2\gamma+\theta)}\lambda_{j+1}^{2\gamma+\theta}a_{j+1}^3;
\end{split}
\end{equation}
\begin{equation}\label{basic-ineq6}
\begin{split}
&\lambda_j^{2\gamma}\lambda_{j+1}^\theta b_ja_{j+1}a_{j+2}\\
=&\ \lambda_j^{2\gamma}\lambda_{j+1}^\theta\left(b_ja_{j+1}^{\frac12}\right)\left(a_{j+1}^{\frac12}a_{j+2}^{\frac14}\right)\left(a_{j+2}^{\frac34}\right)\\
\leq& \frac12 \lambda^{\theta}\lambda_j^{2\gamma+\theta}b_j^2a_{j+1}+\frac14 \lambda^{-2\gamma}\lambda_{j+1}^{2\gamma+\theta}a_{j+1}^2a_{j+2}+\frac14 \lambda^{-4\gamma-\theta}\lambda_{j+2}^{2\gamma+\theta}a_{j+2}^3;
\end{split}
\end{equation}
\begin{equation}\label{basic-ineq7}
\begin{split}
&\lambda_j^{2\gamma}\lambda_{j+1}^\theta b_jb_{j+1}b_{j+2}\\
=&\ \lambda_j^{2\gamma}\lambda_{j+1}^\theta \left(b_jb_{j+1}^{\frac12}\right)  \left(b_{j+1}^{\frac12}b_{j+2}^{\frac14}\right)  \left(b_{j+2}^{\frac34}\right)\\
\leq &\ \frac12\lambda^{\theta}\lambda_j^{2\gamma+\theta} b_j^2b_{j+1}
+\frac14 \lambda^{-2\gamma}\lambda_{j+1}^{2\gamma+\theta} b_{j+1}^2b_{j+2}+\frac14 \lambda^{-4\gamma-\theta}\lambda_{j+2}^{2\gamma+\theta} b_{j+2}^3;
\end{split}
\end{equation}
\begin{equation}\label{basic-ineq8}
\begin{split}
2\lambda_j^{2\gamma+\theta} a_jb_ja_{j+1}
=&\ 2\lambda_j^{2\gamma+\theta} (a_ja_{j+1}^{\frac12})(b_ja_{j+1}^{\frac12})\\
\leq &\ \lambda_j^{2\gamma+\theta}a_j^2a_{j+1} + \lambda_j^{2\gamma+\theta}b_j^2a_{j+1}.
\end{split}
\end{equation}


Applying (\ref{basic-ineq1}), (\ref{basic-ineq2}), (\ref{basic-ineq3}) and (\ref{basic-ineq4-1}) to (\ref{aa-1}), multiplying the constant $c_i$, and adding the shells for $j\geq 1$, we obtain
\begin{equation}\label{aa-2}
\begin{split}
&\frac{d}{dt}\left(c_1\sum_{j=1}^\infty\lambda_j^{2\gamma}a_ja_{j+1}\right)\\
\geq & -\nu c_1(1+\lambda^2)\sum_{j=1}^\infty
\lambda_j^{2\gamma+2}a_ja_{j+1}\\
&+c_1\left(1-\frac12\lambda^{-\frac12(2\gamma+\theta)}-\frac14\lambda^{-4\gamma-\theta}-\frac13\lambda^{-2\gamma}\right)\sum_{j=1}^\infty \lambda_{j}^{2\gamma+\theta}a_{j}^3\\
&-\frac{2}{3}c_1\lambda^{-2\gamma}\sum_{j=1}^\infty\lambda_j^{2\gamma+\theta}b_j^3
+c_1\left(1-\frac12\lambda^{-\frac12(2\gamma+\theta)}\right)\sum_{j=1}^\infty\lambda_j^{2\gamma+\theta}a_jb_j^2\\
&-c_1\left(\frac12\lambda^{-\frac12(2\gamma+\theta)}+\frac12\lambda^\theta+\frac14\lambda^{-2\gamma}\right)\sum_{j=1}^\infty \lambda_j^{2\gamma+\theta}a_j^2a_{j+1}\\
&-\frac12c_1\lambda^{-\frac12(2\gamma+\theta)}\sum_{j=1}^\infty\lambda_j^{2\gamma+\theta}b_j^2a_{j+1}.
\end{split}
\end{equation}
Similarly putting (\ref{ab-1}) together with (\ref{basic-ineq5}), (\ref{basic-ineq6}) and (\ref{basic-ineq7}) gives rise to
\begin{equation}\label{ab-2}
\begin{split}
&\frac{d}{dt}\left(c_2\sum_{j=1}^\infty\lambda_j^{2\gamma}b_ja_{j+1}\right)\\
\geq & -(\mu+\nu\lambda^2)c_2\sum_{j=1}^\infty \lambda_j^{2\gamma+2}b_ja_{j+1}+c_2\sum_{j=1}^\infty\lambda_j^{2\gamma+\theta}a_j^2b_j\\
&+c_2\left(1-\frac14\lambda^{-4\gamma-\theta}\right)\sum_{j=1}^\infty\lambda_j^{2\gamma+\theta}b_j^3\\
&-c_2\left( \frac12\lambda^{-\frac12(2\gamma+\theta)}+\frac14\lambda^{-4\gamma-\theta}\right)\sum_{j=1}^\infty\lambda_j^{2\gamma+\theta}a_j^3\\
&-\frac12c_2\left(\lambda^{-\frac12(2\gamma+\theta)}+\lambda^\theta\right)\sum_{j=1}^\infty\lambda_j^{2\gamma+\theta}b_j^2a_{j+1}-\frac14c_2\lambda^{-2\gamma}\sum_{j=1}^\infty\lambda_j^{2\gamma+\theta}a_j^2a_{j+1}\\
&-c_2\left(\frac12\lambda^\theta+\frac14\lambda^{-2\gamma}\right)\sum_{j=1}^\infty\lambda_j^{2\gamma+\theta}b_j^2b_{j+1}.
\end{split}
\end{equation}
In the end, (\ref{ab-j}) along with (\ref{basic-ineq8}) implies
\begin{equation}\label{ab-j-2}
\begin{split}
&\frac{d}{dt}\left(c_3\sum_{j\geq 1}\lambda_j^{2\gamma}a_jb_j\right)\\
\geq&-(\nu+\mu)c_3\sum_{j\geq 1}\lambda_j^{2\gamma+2}a_jb_j+c_3\left(\lambda^{2\gamma}+1\right)\sum_{j\geq 1}\lambda_j^{2\gamma+\theta}a_j^2b_{j+1}\\
&+c_3\left(\lambda^{2\gamma}-1\right)\sum_{j\geq 1}\lambda_j^{2\gamma+\theta}b_j^2b_{j+1}\\
&-c_3 \sum_{j=1}^\infty \lambda_j^{2\gamma+\theta}a_j^2a_{j+1}-c_3 \sum_{j=1}^\infty \lambda_j^{2\gamma+\theta}b_j^2a_{j+1}.
\end{split}
\end{equation}
Comparing the coefficients of $\sum_{j=1}^\infty \lambda_j^{2\gamma+\theta}a_j^3$, $\sum_{j=1}^\infty \lambda_j^{2\gamma+\theta}b_j^3$, $\sum_{j=1}^\infty \lambda_j^{2\gamma+\theta}a_jb_j^2$, \\
$\sum_{j=1}^\infty \lambda_j^{2\gamma+\theta}a_j^2a_{j+1}$ ,
$\sum_{j=1}^\infty \lambda_j^{2\gamma+\theta}b_j^2a_{j+1}$ and $\sum_{j=1}^\infty \lambda_j^{2\gamma+\theta}b_j^2b_{j+1}$  
on the right hand side of (\ref{a2b2}) and (\ref{aa-2})-(\ref{ab-j-2}), we impose the following conditions for a constant $c_4>0$ 
\begin{equation}\label{para-11}
\begin{split}
&c_1\left(1-\frac12\lambda^{-\frac12(2\gamma+\theta)}-\frac14\lambda^{-4\gamma-\theta}-\frac13\lambda^{-2\gamma}\right)\\
&-c_2\left(\frac12\lambda^{-\frac12(2\gamma+\theta)}+\frac14\lambda^{-4\gamma-\theta}\right)\geq c_4,
\end{split}
\end{equation}
\begin{equation}\label{para-12}
c_2\left(1-\frac14\lambda^{-4\gamma-\theta}\right)-\frac23c_1\lambda^{-2\gamma}\geq  c_4,\\
\end{equation}
\begin{equation}\label{para-13}
c_1\left(1-\frac12\lambda^{-\frac12(2\gamma+\theta)}\right)\geq 0,
\end{equation}
\begin{equation}\label{para-14}
\begin{split}
&2(\lambda^{2\gamma}-1)-c_1\left(\frac12\lambda^{-\frac12(2\gamma+\theta)}+\frac12\lambda^\theta+\frac14\lambda^{-2\gamma}\right)
-\frac14c_2\lambda^{-2\gamma}-c_3\geq 0,
\end{split}
\end{equation}
\begin{equation}\label{para-15}
\begin{split}
&2(\lambda^{2\gamma}-1)-\frac12c_1\lambda^{-\frac12(2\gamma+\theta)}
-\frac12 c_2\left(\lambda^{-\frac12(2\gamma+\theta)}+\lambda^\theta\right)-c_3\geq 0,
\end{split}
\end{equation}
\begin{equation}\label{para-16}
c_3 (\lambda^{2\gamma}-1)-c_2\left(\frac12\lambda^\theta+\frac14\lambda^{-2\gamma}\right)\geq 0.
\end{equation}
We can choose 
$0<c_1=c_2\ll c_3\ll 1$, such that there exists a constant $c_4>0$ with the conditions (\ref{para-11})-(\ref{para-16}) satisfied for $\theta>3$, $\lambda\geq2$ and any $\gamma\in(0,3-\theta)$.  Indeed, we observe that: condition (\ref{para-13}) is automatically satisfied; (\ref{para-14}) and (\ref{para-15}) are satisfied provided 
\[c_1\ll c_3, \ \ \ c_3\leq \frac{\lambda^{2\gamma} -1}{\lambda^\theta+\lambda^{-2\gamma}};\]
while (\ref{para-16}) is satisfied if  
\[c_2\leq \frac{4c_3(\lambda^{2\gamma} -1)}{2\lambda^\theta+\lambda^{-2\gamma}};\]
in the end, we can choose $c_1=c_2$ and 
\[c_4=\min\left\{ c_1\left(1-\lambda^{-\frac12(2\gamma+\theta)}-\frac12 \lambda^{-4\gamma-\theta}-\frac13 \lambda^{-2\gamma}\right), c_1\left(1-\frac14 \lambda^{-4\gamma-\theta}-\frac23 \lambda^{-2\gamma}\right)\right\}\]
which makes (\ref{para-11}) and (\ref{para-12}) valid.

For the constants $c_1,c_2, c_3$ and $c_4$ chosen above, we add (\ref{a2b2}) and (\ref{aa-2})-(\ref{ab-j-2}) to infer
\begin{equation}\label{ab-3}
\begin{split}
\frac{d}{dt}\mathcal L(t)\geq &-\nu(1+\lambda^2)\sum_{j=1}^\infty \lambda_j^{2\gamma+2} a_ja_{j+1}-\mu(1+\lambda^2)\sum_{j=1}^\infty \lambda_j^{2\gamma+2} b_jb_{j+1}\\
&-c_3(\nu+\mu)\sum_{j\geq 1}\lambda_j^{2\gamma+2}a_jb_j-2\nu\|a\|_{\gamma+1}^2-2\mu\|b\|_{\gamma+1}^2\\
&+c_4\sum_{j=1}^\infty \lambda_j^{2\gamma+\theta}a_j^3
+c_4\sum_{j=1}^\infty \lambda_j^{2\gamma+\theta}b_j^3.
\end{split}
\end{equation}
In view of the inequalities in Lemma \ref{le-triple} (i) and (iii) and (\ref{ab-3}), we have
\begin{equation}\label{ab-4}
\begin{split}
\frac{d}{dt}\mathcal L(t)\geq &\left(-2\nu-\nu(1+\lambda^2)\lambda^{-\gamma-1}-\frac12c_3(\nu+\mu)\right)\|a\|_{\gamma+1}^2\\
&+\left(-2\mu-\mu(1+\lambda^2)\lambda^{-\gamma-1}-\frac12c_3(\nu+\mu)\right)\|b\|_{\gamma+1}^2\\
&+c_0c_4\|a\|_{\gamma+1}^3+c_0c_4\|b\|_{\gamma+1}^3\\
\geq &-M_1\left(\|a\|_{\gamma+1}^2+\|b\|_{\gamma+1}^2\right)+\frac12c_0c_4\left(\|a\|_{\gamma+1}^2+\|b\|_{\gamma+1}^2\right)^{\frac32}\\
=& \left(\|a\|_{\gamma+1}^2+\|b\|_{\gamma+1}^2\right)\left(\frac12c_0c_4\left(\|a\|_{\gamma+1}^2+\|b\|_{\gamma+1}^2\right)^{\frac12}-M_1\right)
\end{split}
\end{equation}
where we denote $M_1:=2(\nu+\mu)+(\nu+\mu)(1+\lambda^2)\lambda^{-\gamma-1}+c_3(\nu+\mu)$.  

In the following, we will show that for an appropriate constant $M_0>0$ the assumption $\|a(0)\|_{\gamma}^2+\|b(0)\|_{\gamma}^2>M_0^2$ can close the argument.
Indeed, we define 
\begin{equation}\label{m0}
M_0:=\frac{4M_1}{c_0c_4}(2+2\lambda^{-\gamma-1})^{\frac12}>\frac{4M_1}{c_0c_4}.
\end{equation}
Thus, it follows from the assumption $\|a(0)\|_{\gamma}^2+\|b(0)\|_{\gamma}^2>M_0^2$ that 
\[\|a(0)\|_{\gamma+1}^2+\|b(0)\|_{\gamma+1}^2\geq \|a(0)\|_{\gamma}^2+\|b(0)\|_{\gamma}^2>M_0^2\]
and hence by (\ref{m0}) we have
\[\frac12c_0c_4\left(\|a(0)\|_{\gamma+1}^2+\|b(0)\|_{\gamma+1}^2\right)^{\frac12}-M_1>\frac12c_0c_4M_0-M_1\geq M_1>0.\]
Therefore, (\ref{ab-4}) implies that 
\[\left.\frac{d}{dt}\mathcal L(t)\right\vert_{t=0}>0,\]
and hence, there exists a small time $T>0$ such that
\begin{equation}\label{est-L}
\mathcal L(t)>\mathcal L(0), \ \ \forall t\in(0,T].
\end{equation}

We are left to show that $\mathcal L$ satisfies a Riccati type of inequality. Based on (\ref{ab-4}), we just need to show that 
\begin{equation}\label{m0m1}
\frac14c_0c_4\left(\|a(t)\|_{\gamma+1}^2+\|b(t)\|_{\gamma+1}^2\right)^{\frac12}-M_1\geq 0
\end{equation}
for $t$ on a small time interval, which is expected due to the largeness assumption on the initial data. 
In fact,  from the definition of $\mathcal L(t)$ in (\ref{L1}), it is obviously true for any $t\geq 0$
\begin{equation}\label{Lab}
\|a(t)\|_{\gamma}^2+\|b(t)\|_{\gamma}^2\leq \mathcal L(t)\leq \left(1+(c_1+c_2)\lambda^{-\gamma-1}+c_3\right)\left(\|a(t)\|_{\gamma}^2+\|b(t)\|_{\gamma}^2\right)
\end{equation}
in view of the inequalities of Lemma \ref{le-triple} (iii) with $\gamma+1$ replaced by $\gamma$.
Thus, (\ref{est-L}) and (\ref{Lab}) imply
\begin{equation}\label{Lab1}
\mathcal L(t)\geq \mathcal L(0)\geq \|a(0)\|_{\gamma}^2+\|b(0)\|_{\gamma}^2>M_0^2, \ \ \ \forall \ \ t\in[0,T].
\end{equation}
The right hand side of (\ref{Lab}) also indicates for any $t\geq0$
\begin{equation}\label{Lab2}
\mathcal L(t)\leq \left(2+2\lambda^{-\gamma-1}\right)\left(\|a(t)\|_{\gamma+1}^2+\|b(t)\|_{\gamma+1}^2\right).
\end{equation}
We conclude from (\ref{Lab1}) and (\ref{Lab2})
\[\|a(t)\|_{\gamma+1}^2+\|b(t)\|_{\gamma+1}^2\geq \frac{\mathcal L(t)}{2+2\lambda^{-\gamma-1}}>\frac{M_0^2}{2+2\lambda^{-\gamma-1}}, \ \ t\in[0,T],\]
and hence the definition of $M_0$ in (\ref{m0}) implies (\ref{m0m1}). As a consequence, we have on $[0,T]$ 
\begin{equation}\label{ab-5}
\begin{split}
&\left(\|a(t)\|_{\gamma+1}^2+\|b(t)\|_{\gamma+1}^2\right)\left(\frac12c_0c_4\left(\|a(t)\|_{\gamma+1}^2+\|b(t)\|_{\gamma+1}^2\right)^{\frac12}-M_1\right)\\
\geq&\frac14 c_0c_4\left(\|a(t)\|_{\gamma+1}^2+\|b(t)\|_{\gamma+1}^2\right)^{\frac32}.
\end{split}
\end{equation}
It follows from (\ref{ab-4}), (\ref{ab-5})  and (\ref{Lab2}) that
\begin{equation}\label{LL}
\begin{split}
\frac{d}{dt}\mathcal L(t)\geq & \frac14 c_0c_4\left(\|a(t)\|_{\gamma+1}^2+\|b(t)\|_{\gamma+1}^2\right)^{\frac32}\\
\geq& \frac14 c_0c_4 (2+2\lambda^{-\gamma-1})^{-\frac32}\mathcal L^{\frac32}(t), \ \ t\in[0,T].
\end{split}
\end{equation}
In the end, we point out that since $\mathcal L(T)\geq \mathcal L(0)>M_0^2$, we can start at time $T$ and repeat the process above iteratively to show that 
the Riccati type inequality (\ref{LL}) holds for all $t\geq 0$. It indicates that $\mathcal L(t)$ blows up in finite time. 

\cbdu

\section{Blow-up of positive solutions of dyadic Hall-MHD}
\label{sec-blow-hmhd}


In this section, we prove the blow up of positive solution to the Hall MHD dyadic model (\ref{sys-2}) with $d_i>0$. The strategy of the proof is similar to that of Theorem \ref{thm-blow-mhd} for the MHD dyadic model. 
With the preparation of Lemma \ref{le-cont2}, in order to prove Theorem \ref{thm-blow-hmhd}, it is sufficient to show the following lemma.

\begin{Lemma}\label{claim2}
Consider system (\ref{sys-2}) with $d_i>0$. Let $\theta>3+\gamma$ and $0<\gamma\ll 1$.
Assume $\|a(0)\|_\gamma^2+\|b(0)\|_\gamma^2> M_0^2$ for a certain constant $M_0>0$. The function $\mathcal L(t)$ defined in (\ref{L2}) for positive solution $(a(t), b(t))$ of (\ref{sys-2}) is a Lyapunov function and it blows up in finite time. 
\end{Lemma}
\pf
The main step is to establish a Riccati type inequality for $\mathcal L$. To do so, direct computation based on (\ref{sys-2}) with $d_i>0$ ($d_i=1$ is taken to reduce the number of parameters) gives us
\begin{equation}\label{aa-1h}
\begin{split}
\frac{d}{dt}\left(\lambda_j^{2\gamma}a_ja_{j+1}\right)=&-\nu(1+\lambda^2)\lambda_j^{2\gamma+2}a_ja_{j+1}+\lambda_j^{2\gamma+\theta} a_j^3+\lambda_j^{2\gamma+\theta} a_jb_j^2\\
&+\lambda_{j-1}^\theta \lambda_j^{2\gamma} a_{j-1}^2a_{j+1}+\lambda_{j-1}^\theta \lambda_j^{2\gamma}b_{j-1}^2a_{j+1}\\
&-\lambda_j^{2\gamma+\theta} b_j a_{j+1}b_{j+1}-\lambda_j^{2\gamma}\lambda_{j+1}^\theta a_jb_{j+1}b_{j+2}\\
&-\lambda_j^{2\gamma+\theta} a_ja_{j+1}^2
-\lambda_j^{2\gamma}\lambda_{j+1}^\theta a_ja_{j+1}a_{j+2}, \\
\end{split}
\end{equation}
\begin{equation}\label{bb-1h}
\begin{split}
\frac{d}{dt}\left(\lambda_j^{2\gamma}b_jb_{j+1}\right)=&-\mu(1+\lambda^2)\lambda_j^{2\gamma+2}b_jb_{j+1}+\lambda_j^{2\gamma+\theta+1}b_j^3+\lambda_j^{2\gamma+\theta} a_jb_{j+1}^2
\\
&+\lambda_j^{2\gamma}\lambda_{j+1}^\theta  b_ja_{j+1}b_{j+2}+\lambda_{j-1}^{\theta+1}\lambda_j^{2\gamma}b_{j-1}^2b_{j+1}\\
&-\lambda_{j}^{2\gamma+\theta} b_{j}a_{j+1}b_{j+1}
-\lambda_j^{2\gamma}\lambda_{j+1}^\theta b_j b_{j+1}a_{j+2}\\
&-\lambda_j^{2\gamma+\theta+1}b_jb_{j+1}^2-\lambda_j^{2\gamma}\lambda_{j+1}^{\theta+1}b_jb_{j+1}b_{j+2},
\end{split}
\end{equation}
\begin{equation}\label{ab-1h}
\begin{split}
&\frac{d}{dt}\left(\|a(t)\|_\gamma^2+\|b(t)\|_\gamma^2\right)\\
=&\ -2\nu\|a(t)\|_{\gamma+1}^2-2\mu \|b(t)\|_{\gamma+1}^2+2(\lambda^{2\gamma}-1)\sum_{j=1}^\infty \lambda_j^{2\gamma+\theta}a_j^2a_{j+1}\\
&+2(\lambda^{2\gamma}-1)\sum_{j=1}^\infty \lambda_j^{2\gamma+\theta} b_{j}^2a_{j+1}+2(\lambda^{2\gamma}-1)\sum_{j=1}^\infty\lambda_j^{2\gamma+\theta+1}b_j^2b_{j+1}.\\
\end{split}
\end{equation}
The negative terms on the right hand side of (\ref{aa-1h})-(\ref{ab-1h}) are estimated below, by Young's inequality
\begin{equation}\label{basic-ineq1-h}
\begin{split}
&\lambda_j^{2\gamma+\theta}a_ja_{j+1}^2\\
=&\ \lambda^{-\frac12(2\gamma+\theta)}\left(\lambda_j^{\frac12(2\gamma+\theta)}a_ja_{j+1}^{\frac12}\right)\left(\lambda_{j+1}^{\frac12(2\gamma+\theta)}a_{j+1}^{\frac32}\right)\\
\leq &\ \frac12\lambda^{-\frac12(2\gamma+\theta)} \lambda_{j+1}^{2\gamma+\theta}a_{j+1}^3+\frac12\lambda^{-\frac12(2\gamma+\theta)} \lambda_{j}^{2\gamma+\theta}a_j^2a_{j+1};
\end{split}
\end{equation}
\begin{equation}\label{basic-ineq2-h}
\begin{split}
&\lambda_j^{2\gamma+\theta}b_ja_{j+1}b_{j+1}\\
=&\ \lambda^{-\frac12(2\gamma+\theta)} \left(\lambda_j^{\frac12(2\gamma+\theta)}b_ja_{j+1}^{\frac12}\right) \left(\lambda_{j+1}^{\frac12(2\gamma+\theta)}a_{j+1}^{\frac12}b_{j+1}\right)\\
\leq&\  \frac12 \lambda^{-\frac12(2\gamma+\theta)} \lambda_j^{2\gamma+\theta}b_j^2a_{j+1}+
\frac12 \lambda^{-\frac12(2\gamma+\theta)} \lambda_{j+1}^{2\gamma+\theta}a_{j+1}b_{j+1}^2;
\end{split}
\end{equation}
\begin{equation}\label{basic-ineq3-h}
\begin{split}
&\lambda_j^{2\gamma}\lambda_{j+1}^\theta a_ja_{j+1}a_{j+2}\\
=&\ \lambda_j^{2\gamma}\lambda_{j+1}^\theta \left( a_ja_{j+1}^{\frac12}\right)
\left(a_{j+1}^{\frac12}a_{j+2}^{\frac14}\right)\left(a_{j+2}^{\frac34}\right)\\
\leq &\ \frac12\lambda^\theta \lambda_j^{2\gamma+\theta}a_j^2a_{j+1} +\frac14\lambda^{-2\gamma} \lambda_{j+1}^{2\gamma+\theta}a_{j+1}^2a_{j+2}+\frac14 \lambda^{-4\gamma-\theta}\lambda_{j+2}^{2\gamma+\theta}a_{j+2}^3;
\end{split}
\end{equation}
\begin{equation}\label{basic-ineq4-h}
\begin{split}
&\lambda_j^{2\gamma}\lambda_{j+1}^\theta a_jb_{j+1}b_{j+2}\\
=&\ \lambda^{-(2\gamma+1)}\lambda_j^{-\frac23}\left(\lambda_j^{\frac13(2\gamma+\theta)}a_j\right) \left(\lambda_{j+1}^{\frac13(2\gamma+\theta+1)}b_{j+1}\right) \left(\lambda_{j+2}^{\frac13(2\gamma+\theta+1)}b_{j+2}\right)\\
\leq &\ \frac13 \lambda^{-(2\gamma+\frac53)} \lambda_j^{2\gamma+\theta}a_j^3+\frac13 \lambda^{-(2\gamma+\frac53)} \lambda_{j+1}^{2\gamma+\theta+1}b_{j+1}^3\\
&+\frac13 \lambda^{-(2\gamma+\frac53)} \lambda_{j+2}^{2\gamma+\theta+1}b_{j+2}^3;
\end{split}
\end{equation}
\begin{equation}\label{basic-ineq5-h}
\begin{split}
&\lambda_j^{2\gamma}\lambda_{j+1}^\theta b_jb_{j+1}a_{j+2}\\
=&\ \lambda_j^{2\gamma}\lambda_{j+1}^{\theta}\left(b_jb_{j+1}^{\frac12}\right) \left(b_{j+1}^{\frac12}\right) \left(a_{j+2}\right)\\
\leq &\ \frac12 \lambda_j^{2\gamma}\lambda_{j+1}^{\theta} b_j^2b_{j+1}+\frac16\lambda_j^{2\gamma}\lambda_{j+1}^{\theta} b_{j+1}^3+\frac13 \lambda_j^{2\gamma}\lambda_{j+1}^{\theta} a_{j+2}^3\\
\leq &\ \frac12 \lambda^{\theta-1}\lambda_{j}^{2\gamma+\theta+1} b_j^2b_{j+1}+\frac16\lambda^{-2\gamma-2}\lambda_{j+1}^{2\gamma+\theta+1} b_{j+1}^3\\
&+\frac13 \lambda^{-4\gamma-\theta}\lambda_{j+2}^{2\gamma+\theta} a_{j+2}^3;
\end{split}
\end{equation}
\begin{equation}\label{basic-ineq6-h}
\begin{split}
&\lambda_j^{2\gamma+\theta+1}b_jb_{j+1}^2\\
=&\lambda_j^{2\gamma+\theta+1}\left(b_jb_{j+1}^{\frac12}\right)
\left(b_{j+1}^{\frac32}\right)\\
\leq & \frac12 \lambda_j^{2\gamma+\theta+1}b_j^2b_{j+1}
+ \frac12 \lambda^{-(2\gamma+\theta+1)}\lambda_{j+1}^{2\gamma+\theta+1}b_{j+1}^3;
\end{split}
\end{equation}
\begin{equation}\label{basic-ineq7-h}
\begin{split}
&\lambda_j^{2\gamma}\lambda_{j+1}^{\theta+1}b_jb_{j+1}b_{j+2}\\
=&\lambda_j^{2\gamma}\lambda_{j+1}^{\theta+1}\left(b_jb_{j+1}^{\frac12} \right)
\left(b_{j+1}^{\frac12}b_{j+2}^{\frac14}\right)\left(b_{j+2}^{\frac34}\right)\\
\leq&\ \frac12 \lambda^{\theta+1}\lambda_j^{2\gamma+\theta+1} b_j^2b_{j+1}+\frac14\lambda^{-2\gamma}\lambda_{j+1}^{2\gamma+\theta+1} b_{j+1}^2b_{j+2}\\
&+\frac14\lambda^{-4\gamma-\theta-1}\lambda_{j+2}^{2\gamma+\theta+1} b_{j+2}^3.
\end{split}
\end{equation}
Applying (\ref{basic-ineq1-h})-(\ref{basic-ineq4-h}) to (\ref{aa-1h}) yields
\begin{equation}\label{aa-2h}
\begin{split}
&\frac{d}{dt}\left(c_1\sum_{j=1}^\infty\lambda_j^{2\gamma}a_ja_{j+1}\right)\\
\geq & -c_1\nu(1+\lambda^2)\sum_{j=1}^\infty
\lambda_j^{2\gamma+2}a_ja_{j+1}\\
&+c_1\left(1-\frac13\lambda^{-(2\gamma+\frac53)}-\frac12\lambda^{-\frac12(2\gamma+\theta)}-\frac14\lambda^{-4\gamma-\theta}\right)\sum_{j=1}^\infty \lambda_{j}^{2\gamma+\theta}a_{j}^3\\
&+c_1\left(1-\frac12\lambda^{-\frac12(2\gamma+\theta)}\right) \sum_{j=1}^\infty \lambda_{j}^{2\gamma+\theta}a_jb_j^2\\
&-\frac{2}{3}c_1 \lambda^{-(2\gamma+\frac53)}\sum_{j=1}^\infty\lambda_j^{2\gamma+\theta+1}b_j^3
-\frac12c_1\lambda^{-\frac12(2\gamma+\theta)} \sum_{j=1}^\infty \lambda_j^{2\gamma+\theta}b_j^2a_{j+1}\\
&-c_1\left(\frac12\lambda^{-\frac12(2\gamma+\theta)}+\frac12\lambda^\theta+\frac14\lambda^{-2\gamma}\right)\sum_{j=1}^\infty \lambda_j^{2\gamma+\theta}a_j^2a_{j+1}.
\end{split}
\end{equation}
While (\ref{basic-ineq2-h}) and (\ref{basic-ineq5-h})-(\ref{basic-ineq7-h}) applied to (\ref{bb-1h}) gives
\begin{equation}\label{bb-2h}
\begin{split}
&\frac{d}{dt}\left(c_2\sum_{j=1}^\infty\lambda_j^{2\gamma}b_jb_{j+1}\right)\\
\geq & -c_2\mu(1+\lambda^2)\sum_{j=1}^\infty
\lambda_j^{2\gamma+2}b_jb_{j+1}\\
&-\frac13c_2\lambda^{-4\gamma-\theta}\sum_{j=1}^\infty \lambda_{j}^{2\gamma+\theta}a_{j}^3
-\frac12c_2\lambda^{-\frac12(2\gamma+\theta)} \sum_{j=1}^\infty \lambda_{j}^{2\gamma+\theta}a_jb_j^2\\
&+c_2 \left(1-\frac16\lambda^{-2\gamma-2}-\frac12\lambda^{-2\gamma-\theta-1}-\frac14\lambda^{-4\gamma-\theta-1}\right)\sum_{j=1}^\infty\lambda_j^{2\gamma+\theta+1}b_j^3\\
&-\frac12c_2\lambda^{-\frac12(2\gamma+\theta)}\sum_{j=1}^\infty \lambda_j^{2\gamma+\theta}b_j^2a_{j+1}\\
&-c_2\left(\frac12\lambda^{\theta-1}+\frac12+\frac12\lambda^{\theta+1}+\frac14\lambda^{-2\gamma} \right) \sum_{j=1}^\infty \lambda_j^{2\gamma+\theta+1}b_j^2b_{j+1}.
\end{split}
\end{equation}
In order to have the negative terms in (\ref{aa-2h})-(\ref{bb-2h}) and (\ref{ab-1h}) absorbed by the positive terms, we claim there exists a constant $c_3>0$ such that
\begin{equation}\label{para-11-h}
\begin{split}
&c_1\left(1-\frac13\lambda^{-(2\gamma+\frac53)}-\frac12\lambda^{-\frac12(2\gamma+\theta)}-\frac14\lambda^{-4\gamma-\theta}\right)-\frac13c_2\lambda^{-4\gamma-\theta}\geq c_3,\\
\end{split}
\end{equation}
\begin{equation}\label{para-12-h}
\begin{split}
&c_2\left(1-\frac16\lambda^{-2\gamma-2}-\frac12\lambda^{-2\gamma-\theta-1}-\frac14\lambda^{-4\gamma-\theta-1}\right)-\frac23c_1\lambda^{-2\gamma-\frac53}\geq c_3,\\
\end{split}
\end{equation}
\begin{equation}\label{para-13-h}
c_1\left(1-\frac12\lambda^{-\frac12(2\gamma+\theta)}\right)-\frac12c_2\lambda^{-\frac12(2\gamma+\theta)}\geq 0,
\end{equation}
\begin{equation}\label{para-14-h}
2(\lambda^{2\gamma}-1)-c_1\left(\frac12\lambda^{-\frac12(2\gamma+\theta)}+\frac12\lambda^\theta+\frac14\lambda^{-2\gamma} \right)\geq 0,
\end{equation}
\begin{equation}\label{para-15-h}
2(\lambda^{2\gamma}-1)-\frac12 c_1\lambda^{-\frac12(2\gamma+\theta)} -\frac12 c_2\lambda^{-\frac12(2\gamma+\theta)}\geq 0,
\end{equation}
\begin{equation}\label{para-16-h}
2(\lambda^{2\gamma}-1)-c_2\left(\frac12\lambda^{\theta-1}+\frac12+\frac12\lambda^{\theta+1}+\frac14\lambda^{-2\gamma} \right)\geq 0.
\end{equation}
As a matter of fact, we can choose $c_2=\frac12c_1$ and $0<c_1\ll 1$ such that 
\begin{equation}\label{para-c1-h}
c_1\leq \frac{8(\lambda^{2\gamma}-1)}{2\lambda^{\theta-1}+2+2\lambda^{\theta+1}+\lambda^{-2\gamma} }.
\end{equation}
One can check conditions (\ref{para-13-h})-(\ref{para-16-h}) are satisfied. Consequently, for $\lambda\geq 2$, there exists a constant $c_3>0$ such that  (\ref{para-11-h}) and (\ref{para-12-h}) are also satisfied.

In view of (\ref{L2}), adding (\ref{ab-1h}) and (\ref{aa-2h})-(\ref{bb-2h}) leads to
\begin{equation}\label{ab-3h}
\begin{split}
\frac{d}{dt}\mathcal L(t)\geq &-c_1\nu(1+\lambda^2)\sum_{j=1}^\infty \lambda_j^{2\gamma+2} a_ja_{j+1}-c_2\mu(1+\lambda^2)\sum_{j=1}^\infty \lambda_j^{2\gamma+2} b_jb_{j+1}\\
&-2\nu\|a\|_{\gamma+1}^2-2\mu\|b\|_{\gamma+1}^2+c_3\sum_{j=1}^\infty \lambda_j^{2\gamma+\theta}a_j^3
+c_3\sum_{j=1}^\infty \lambda_j^{2\gamma+\theta+1}b_j^3.
\end{split}
\end{equation}
Applying the inequalities of Lemma \ref{le-triple} to (\ref{ab-3h}), we obtain
\begin{equation}\label{ab-4h}
\begin{split}
\frac{d}{dt}\mathcal L(t)\geq &\left(-2\nu-c_1\nu(1+\lambda^2)\lambda^{-\gamma-1}\right)\|a\|_{\gamma+1}^2\\
&+\left(-2\mu-c_2\mu(1+\lambda^2)\lambda^{-\gamma-1}\right)\|b\|_{\gamma+1}^2\\
&+c_0c_3\|a\|_{\gamma+1}^3+c_0c_3\|b\|_{\gamma+1}^3\\
\geq &-M_1\left(\|a\|_{\gamma+1}^2+\|b\|_{\gamma+1}^2\right)+\frac12c_0c_3\left(\|a\|_{\gamma+1}^2+\|b\|_{\gamma+1}^2\right)^{\frac32}\\
=& \left(\|a\|_{\gamma+1}^2+\|b\|_{\gamma+1}^2\right)\left(\frac12c_0c_3\left(\|a\|_{\gamma+1}^2+\|b\|_{\gamma+1}^2\right)^{\frac12}-M_1\right)
\end{split}
\end{equation}
with $M_1:=2(\nu+\mu)+(c_1\nu+c_2\mu)(1+\lambda^2)\lambda^{-\gamma-1}$.  Define 
\begin{equation}\notag
M_0:=\frac{4M_1}{c_0c_3}(1+(c_1+c_2)\lambda^{-\gamma-1})^{\frac12}>\frac{4M_1}{c_0c_3}.
\end{equation}
With such $M_0$ and the estimate (\ref{ab-4h}), an analogous analysis as the last part of the proof of Lemma \ref{le-blow-mhd} can be used to justify the statement of the current lemma.

\cbdu



\section*{Acknowledgement}
The author is sincerely grateful to Professor Susan Friedlander for sharing some references on dyadic MHD models from the physics community, and for many insightful conversations with her.

\end{document}